\documentclass[10pt]{article}
\usepackage{amsmath,amssymb}
\newtheorem{thm}{Theorem}[section]
\newtheorem{co}[thm]{Corollary}
\newtheorem{lem}[thm]{Lemma}

\newtheorem{definition}[thm]{Definition}

\newtheorem{example}[thm]{Example}
\newenvironment{exmp}{\begin{example}\rm}{\end{example}}
\newtheorem{remark}[thm]{Remark}

\newcommand{\N}{\mathbb{N}}
\newcommand{\R}{\mathbb{R}}
\newcommand{\trace}{\mathrm{trace}\,}

\newcommand{\Section}[1]{\section{#1}\setcounter{equation}{0}}

\newcommand{\openbox}{\leavevmode
  \hbox to.77778em{%
    \hfil\vrule
  \vbox to.675em{\hrule width.6em\vfil\hrule}%
  \vrule\hfil}}
\newcommand{\proofname}{Proof}

\newenvironment{proof}[1][\proofname]{\par\normalfont
  \trivlist\item[\hskip\labelsep\itshape #1:]\ignorespaces
  }{\hspace*{1cm}\hspace*{\fill}\openbox \medskip\endtrivlist}

\title{A Note on the Weighted Harmonic-Geometric-Arithmetic Means
  Inequalities
  \footnote{To appear in Math. Ineq. App.}
}%
\date{\today}%
\author{G\'erard Maze, Urs Wagner\\
{\small {\em e-mail:\/} \{gmaze,uwagner\}@math.uzh.ch \vspace{-1mm} }\\
{\small Mathematics Institute\vspace{-1mm}}\\
{\small University of Z\"urich\vspace{-1mm}}\\
{\small Winterthurerstr 190, CH-8057 Z\"urich, Switzerland }
\vspace{3mm} }

\begin{document}\maketitle
\thispagestyle{empty}
\begin{abstract}
  In this note, we derive non trivial sharp bounds related to the
  weighted harmonic-geometric-arithmetic means inequalities, when two
  out of the three terms are known. As application, we give an
  explicit bound for the trace of the inverse of a symmetric positive
  definite matrix and an inequality related to the coefficients of
  polynomials with positive roots.
\end{abstract}
\vspace{3mm}

\Section{Introduction and Main Results}

The well known weighted harmonic-geometric-arithmetic means
inequalities (HGA) can be stated as follows. Let $\alpha_i>0$ and $x_i >
0$, $i=1,\ldots,n$ with $\sum_i \alpha_i = 1$, and define $h,g,a$ by
\begin{equation*}
h=\left( \sum_{i=1}^n \frac{\alpha_i}{x_i} \right)^{-1},  \hspace{5mm}
g=\prod_{i=1}^n x_i^{\alpha_i}, \hspace{5mm}
a=\sum_{i=1}^n \alpha_i x_i.
\end{equation*}
Then the HGA inequalities state that
\begin{equation}\label{HGA}
h \leq g \leq a.
\end{equation}
One equality is reached if and only if all the $x_i$ are equal, which
then implies that both inequalities are in fact equalities.  The
terms of the previous inequalities are respectively called the
harmonic, the geometric and the arithmetic mean of the $x_i$ with
weight $\alpha_i$. There exist several extensions of these
inequalities, see for example \cite{ineq1,ineq2,ineq,ineq3}. In this
note we focus on the case where two of the means are known and non
trivial bounds on the third have to be determined. Actually, Theorem
\ref{GA} below gives a sharp lower bound and a sharp upper bound on
the harmonic mean, when both the arithmetic and the geometric means
are known. The dual bounds, i.e., an upper and a lower bound on the
arithmetic mean when both the harmonic and the geometric means are
known can easily be deduced with the change of variables
$y_i=x_i^{-1}$. Theorem \ref{HA} gives a sharp lower bound and a sharp
upper bound on the geometric mean, when both the harmonic and the
arithmetic means are known, extending Inequalities (\ref{HGA}) above
when the two extreme values are in fact known.\\

The theory of complementary inequalities is a field, where upper
bounds for the ratios $a/g$, $a/h$, $g/h$ and for the differences
$a-g$, $a-h$, $g-h$ are obtained in terms of the upper and lower
bounds for the variables $x_i$.  For instance, Kantorovich's
inequality, see e.g. \cite{ineq3}, provides a well-known upper bound
for $a/h$ under these conditions.  Despite our search for a similar
result in the vast literature on the subject, we were not able to find
the inequalities presented in this article in any published work.
Related to our results, let us however mention \cite{pales}, where the
author considers the interplay of the three means $h, g, a$. In this
paper, it is shown that the moment space of the triplets $(h, g, a)$
is the set $M = \{(u, v, w) \in \R^3 \; : \; 0 \leq u \leq v \leq
w\}$. This means that for any positive $\varepsilon$, there exists
$n\in \N$, $x_1, \ldots,x_n > 0$ such that
\[
|h-u| \leq \varepsilon, \; |g-v| \leq \varepsilon, \; |a-w| \leq
\varepsilon.
\]
The meaning of the result is that if $n$ is not fixed, then the only
meaningful inequality for the three means is $h \leq g \leq a$. In the
present paper, $n$ is fixed and we will suppose that at least one
$x_i$ is different from the others, which insure that Inequalities
(\ref{HGA}) are strict. The main results are the following.

\begin{thm}\label{GA}
With the above notations, if $\alpha = \min_i\{\alpha_i\}$ and $\xi_0
\in [0,1]$, $\xi_1 \in [1,1/\alpha]$ are the solutions of the
equation
\[
g=a \, \xi^{\alpha}\left( \frac{1-\alpha \xi}{1-\alpha
}\right)^{1-\alpha},
\]
then 
\[
a  \, \frac{\xi_0(1-\alpha \xi_0)}{\alpha-2\alpha\xi_0+\xi_0} 
\leq h \leq
a  \, \frac{\xi_1(1-\alpha \xi_1)}{\alpha-2\alpha\xi_1+\xi_1}.
\]
The first (resp. second) inequality reaches equality if and only if
$x_j=\xi_0$ (resp. $x_j=\xi_1$) and $x_l=x_k$, $\forall l,k\neq j$ for
some $j$ with $\alpha_j= \min_i\{\alpha_i\}$.
\end{thm}
The uniqueness of the solutions $\xi_0$ and $\xi_1$ will be made clear
in the sequel. Based on this result, we give explicit general lower
and upper bounds for the harmonic and arithmetic means in Corollary
\ref{co}. As application, we give an explicit bound for the trace of the
inverse of a symmetric positive definite matrix in Example \ref{ex3}
and for the quotient of coefficients of polynomials with positive
roots in Example \ref{ex4}.

\begin{thm}\label{HA}
With the above notations, if  $\alpha = \min_i\{\alpha_i\}$
and  $\xi_0 \in [0,1]$, $\xi_1 \in [1,1/\alpha]$ are the solutions
of the equation
\[
h=a \, \frac{\xi(1-\alpha \xi)}{\alpha-2\alpha\xi+\xi},
\]
then we have
\[
 a \, \xi_1^{\alpha}\left( \frac{1-\alpha
 \xi_1}{1-\alpha}\right)^{1-\alpha}
 \,\leq  \, g  \,\leq \,
 a \, \xi_0^{\alpha}\left( \frac{1-\alpha
 \xi_0}{1-\alpha}\right)^{1-\alpha}.
\] 
The first (resp. second) inequality reaches equality if and only if
$x_j=\xi_1$ (resp. $x_j=\xi_0$) and $x_l=x_k$, $\forall l,k\neq j$
for some $j$ with $\alpha_j= \min_i\{\alpha_i\}$.
\end{thm}
Based on this result, we give explicit sharp lower and upper bounds
for the geometric mean in Corollary \ref{explicit} and simpler bounds
in Corollary \ref{co}.

\Section{Explicit Bounds}

We postpone the proof of Theorem \ref{GA} and Theorem \ref{HA} and
present explicit bounds for the different means. So let $\alpha_i>0$
with $\sum \alpha_i=1$ and $x_i> 0$, $i=1,\ldots,n$ be real numbers
such that
\begin{equation*}
h=\left( \sum_{i=1}^n \frac{\alpha_i}{x_i} \right)^{-1},  \hspace{5mm}
g=\prod_{i=1}^n x_i^{\alpha_i}, \hspace{5mm}
a=\sum_{i=1}^n \alpha_i x_i
\end{equation*}
and let $\alpha = \min_i\{\alpha_i\}$. Note that $\sum_{i=1}^n
\alpha_i = 1$ implies that $\alpha \leq 1/n$. In the case of Theorem
\ref{HA}, the equation for $\xi$ is exactly solvable, and one readily
verifies that a sharp bound in closed form can be computed as follows.

\begin{co}\label{explicit}
With the above notations, we have $\alpha \leq 1/n$ and 
\begin{eqnarray*}
g & \leq & \left(
\frac{a-h(1-2\alpha)-\sqrt{(a-h)(a-h(1-2\alpha)^2)}}{2\alpha}\right)^{\alpha}
\left(
\frac{a+h(1-2\alpha)+\sqrt{(a-h)(a-h(1-2\alpha)^2)}}{2(1-\alpha)}\right)^{1-\alpha},\\
g & \geq & \left(
\frac{a-h(1-2\alpha)+\sqrt{(a-h)(a-h(1-2\alpha)^2)}}{2\alpha}\right)^{\alpha}
\left(
\frac{a+h(1-2\alpha)-\sqrt{(a-h)(a-h(1-2\alpha)^2)}}{2(1-\alpha)}\right)^{1-\alpha}.
\end{eqnarray*}
\end{co}
The bounds of the next corollary are not sharp anymore but are both in
closed form and simple.
\begin{co}\label{co}
With the above notations, we have
\[
a \cdot \left (\alpha e
\left(\frac{a}{g}\right)^{1/\alpha} + 1 \right)^{-1}
< h \leq g \leq a <
h \cdot \left (\alpha e
\left(\frac{g}{h}\right)^{1/\alpha} + 1 \right),
\]
and
\[
h \cdot \left(\frac{h}{a}\exp \left( \frac{h}{a} +\frac{n}{n-1} \right)
\right)^{-\alpha}
< g <
a \cdot \left(\frac{h}{a}\exp \left( \frac{h}{a} +\frac{n}{n-1} \right)
\right)^{\alpha}.
\]
\end{co}
Asymptotically with $n$, the last two inequalities give an improvement
of the usual HGA inequalities when $h/a < t_0 = 0.278464...$, where
$t_0 e^{t_0+1}=1$.
\begin{proof}
Let us start with the first set of inequalities. The root
$\xi=\xi_0<1$ of Theorem \ref{GA} satisfies the following
inequalities,
\[
(g/a)^{1/\alpha} =\xi\left( \frac{1-\alpha
  \xi}{1-\alpha }\right)^{\frac{1-\alpha}{\alpha}} =\xi\left( 1+\frac{\alpha}
  {1-\alpha } (1-\xi)\right)^{\frac{1-\alpha}{\alpha}}< \xi e^{1-\xi}
  < \xi e
\]
because $\frac{\alpha}{1-\alpha}<1$ and $(1+u)^v<e^{uv}$, as soon as
$v>0$ and $|u|<1$.  Now, since $\frac{1-\alpha\xi}{1-\alpha}>1$ and
$1/\xi< e (a/g)^{1/\alpha}$,
\[
\frac{1}{h}=\sum_{i=1}^n \frac{\alpha_i}{x_i}  =\frac{1}{a} \cdot  
\sum_{i=1}^n 
\frac{\alpha_i}{z_i}  \leq  \frac{1}{a} \cdot \left( 
\frac{\alpha-2\alpha\xi+\xi}{\xi(1-\alpha \xi)} \right) 
=  \frac{1}{a} \cdot \left( \frac{\alpha}{\xi} +
\frac{(1-\alpha)^2}{1-\alpha \xi} \right) <  \frac{1}{a} \cdot
\left( \alpha e
(a/g)^{1/\alpha} + 1 \right) 
\]
which is equivalent to $h> a \cdot \left(\alpha e (a/g)^{1/\alpha} + 1
\right)^{-1}$. By setting $x_i'=1/x_i$, we have $a'=1/h, g'=1/g,
h'=1/a$ and the inequality $a< h \cdot \left( \alpha e \left( g/h
\right)^{1/\alpha} + 1 \right)$ is a direct consequence. Let us now
prove the second set of inequalities. Since $(1-2\alpha)^2\leq 1$, we
have
\[
(a-h)^2 \leq (a-h)(a-h(1-2\alpha)^2) \leq (a-h(1-2\alpha)^2)^2 
\]
an thus the upper bound of Corollary \ref{explicit} gives
\begin{eqnarray*}
g & \leq & \left(
\frac{a-h(1-2\alpha)-(a-h)}{2\alpha}\right)^{\alpha}
\left(
\frac{a+h(1-2\alpha)+(a-h(1-2\alpha)^2)}{2(1-\alpha)}\right)^{1-\alpha}.
\end{eqnarray*}
Since $\alpha\leq 1/n$, $\frac{1}{1-\alpha}\leq
1+\frac{n}{n-1}\alpha$, and after suitable simplifications, using once
again the above exponential inequality, we obtain
\begin{eqnarray*}
g & \leq & 
a \cdot \left(\frac{h}{a}\right)^{\alpha} \cdot
\left(1+\frac{h}{a}(1-2\alpha)\alpha \right)^{1-\alpha} \cdot 
\left(1+\frac{n}{n-1}\alpha\right)^{1-\alpha} 
< a \cdot \left(\frac{h}{a}\exp \left( \frac{h}{a} +\frac{n}{n-1}\right) 
\right)^{\alpha}.
\end{eqnarray*}
The reverse inequality is once again obtained by setting
$z_i=1/x_i$. This finishes the proof of the lemma.
\end{proof}

\Section{The Case $n=2$}

For the rest of the article, without loss of generality, we will
assume that the $x_i$ are normalized so that the arithmetic mean is
equal to $1$. This is simply obtained by the change of variable
$x_i'=x_i/a$, leading to $a'=1, g'=g/a$ and $h'=h/a$. \\

\noindent Along the way of the proofs of the main results, we start
with the case $n=2$. This will turn out to be in fact the most
important case, as the general case will be a consequence of it. The
next two lemmas will be useful in the sequel.
\begin{lem} \label{f-g}
Let $\alpha \in ]0,1/2]$ and $f$ and $\varphi$ be the functions
    defined over $[0,1/\alpha]$ defined by
\begin{equation}\label{fg}
f(x) = x^{\alpha}\left( \frac{1-\alpha x}{1-\alpha
}\right)^{1-\alpha} \;\; \mbox{ and } \;\; 
\varphi(x)=\sqrt{\frac{x(1-\alpha x)}{(1-2\alpha)x+\alpha}}.
\end{equation}
Then $\varphi(0)=\varphi(1/\alpha)=f(0)=f(1/\alpha)=0$,
$f(1)=\varphi(1)=1$, they are strictly increasing over $[0,1]$ and
strictly decreasing over $[1,1/\alpha]$, and fulfill the property that
$f > \varphi$ over $[0,1[$, and $f < \varphi$ over $]1,1/\alpha]$.
\end{lem}
\begin{proof}
  Clearly, $\varphi(0)=\varphi(1/\alpha)=f(0)=f(1/\alpha)=0$,
  $g(1)=f(1)=1$. A short analysis of $f$ and of the radical of
  $\varphi$ shows that they are strictly increasing over $[0,1]$ and
  strictly decreasing over $[1,1/\alpha]$. The less obvious fact is
  that $f > \varphi$ over $[0,1[$ and $f < \varphi$ over
      $]1,1/\alpha]$. In order to prove it, let us consider the
  quotient $f/\varphi$. Since $(f/\varphi)(1)=1$, the statement would
  be proved if we can show that $f/\varphi$ is strictly decreasing
  over $[0,1/\alpha]$. Let us prove that it is the case by showing
  that $(f/\varphi)' < 0$ over $]0,1[ \cup ]1,1/\alpha[ $. First
  \[ 
  \left( f/\varphi \right) (x) = \frac{x^{\alpha}\left( \frac{1-\alpha
      x}{1-\alpha }\right)^{1-\alpha}}{\frac{x^{1/2}(1-\alpha
      x)^{1/2}}{((1-2\alpha)x+\alpha)^{1/2}}} = (1-\alpha)^{\alpha-1}
  \cdot \left(\frac{1}{x}-\alpha \right)^{1/2-\alpha} \cdot
  \left((1-2\alpha)x+\alpha \right)^{1/2}.
  \]
  After suitable simplifications, we obtain
  \[
  \left(f/\varphi\right)'(x) =-\alpha(1-\alpha)^{\alpha-1} \cdot
  \frac{(x-1)^2}
  {2x^2(\frac{1}{x}-\alpha)^{\alpha+1/2}(\alpha+(1-2\alpha)x)^{1/2}}.
  \] 
  Note that $\frac{1}{x}-\alpha > 0$ and $1-2\alpha \geq 0$ so the
  condition $(f/\varphi)' < 0$ is fulfilled.
\end{proof}

\begin{lem} \label{gamma}
If $\alpha\in ]0,1[$ and $x\in [0,1/\alpha]$, then
\begin{equation}
\frac{1-x}{1-\alpha x}+
 \ln \left( 1-\frac{1-x}{1-\alpha x}\right) +
\frac{(1-x)^2}{(1-\alpha x)(1-(2\alpha-1)x)} 
\;\;\left\{
\begin{array}{rl}
\,\leq 0 & \mbox{ if } x \in [0,1],\\
\,\geq 0 & \mbox{ if } x \in [1,1/\alpha].\\
\end{array}
\right.
\end{equation}
\end{lem}
\begin{proof}
If $t=\frac{1-x}{1-\alpha x}$, then $-\infty < t \leq 0$ for $x \in
\lbrack 1,1/\alpha \lbrack$ and $0 \leq t \leq 1$ for $x \in \lbrack
0,1 \rbrack$. Since $\frac{1-(2\alpha-1)x}{1- \alpha x}=2-t$, the
above expression is equal to $s(t)=\frac{2t}{2-t } + \ln \left(
1-t\right)$. But since $s(0)=0$ and $s$ has a non positive derivative
$s'(t)=\frac{-t^2}{(2-t)^2(1-t)} \leq 0$, the function $s$ is
decreasing and $s\leq 0$ for $x \in [0,1]$ and $s \geq 0$ for $x \in
[1,1/\alpha[$.
\end{proof}
Returning to the original problem, let us focus on the case where both
$g$ and $a$ are known and an upper and a lower bound on $h$ is to be
determined.  If $\alpha_1,\alpha_2>0$ and $\alpha_1 + \alpha_2=1$, up to a
permutation of the indices, we can assume without loss of generality
that $\alpha_1\leq 1/2$. The two dimensional case can be stated as
follows: given two real numbers $0<\alpha \leq 1/2$ and $0<g<1$, we
want to find the minimal and the maximal value of
\[
H(x,y) = \left(\alpha /x + (1-\alpha) /y\right)^{-1},
\]
where $x$ and $y$ fulfill the conditions
\begin{equation*}
\alpha x + (1-\alpha) y  =  1 \; \mbox{ and } \;
x^{\alpha}y^{1-\alpha}  =  g.
\end{equation*}
Clearly, these conditions imply that
\begin{equation}\label{f}
f(x)=g\mbox{ where } f(x) = x^{\alpha}\left( \frac{1-\alpha x}{1-\alpha
}\right)^{1-\alpha}.
\end{equation}
Note that the function $f$ appears in Lemma \ref{f-g}. We call $x_1$
and $x_2$ the two unique solutions of Equation (\ref{f}), with
$x_1<1<x_2$. Then, with $\varphi$ being the function of Lemma
\ref{f-g},
\[ 
H(x_i,y_i) = \left( \frac{\alpha}{x_i} + \frac{1-\alpha}{y_i}
\right)^{-1} = \frac{x_i(1-\alpha x_i)}{x_i-2x_i \alpha +\alpha} =
\varphi^2(x_i)
\]
and Lemma \ref{f-g} implies that $\varphi(x_2) > f(x_2) = f(x_1) >
\varphi(x_1)$ because of the respective positions of $f$ and
$\varphi$. This directly gives the following lemma:
\begin{lem} \label{2-case ga}
Let $0<\alpha\leq1/2$ and $0<g<1$. If $x$ and $y$ fulfill the
conditions 
\begin{equation*} 
\alpha x + (1-\alpha) y  =  1 \; \mbox{ and } \;
x^{\alpha}y^{1-\alpha}  =  g,
\end{equation*}
then 
\[ 
\frac{x_1(1-\alpha x_1)}{x_1-2x_1 \alpha +\alpha}\leq \left(
\frac{\alpha}{x} + \frac{1-\alpha}{y} \right)^{-1} \leq
\frac{x_2(1-\alpha x_2)}{x_2-2x_2 \alpha +\alpha}
\]
where $x_1$ and $x_2$ are the unique solutions over $[0,1]$ and
$[1,1/\alpha]$ respectively of the equation
\begin{equation}\label{ux}
g= x^{\alpha}\left( \frac{1-\alpha x}{1-\alpha
}\right)^{1-\alpha}.
\end{equation}
\end{lem}
We would like now to prove that for a fixed $g$, and as a function of
$\alpha \in ]0,1/2]$, the minimum and the maximum values above
$H(x_1)$ and $H(x_2)$ are increasing and decreasing functions
respectively. This result will be useful in the sequel. More
precisely, if we set
\[
M(x,\alpha) = \left( \frac{\alpha}{x}
+ \frac{(1-\alpha)^2}{1-\alpha x} \right)^{-1}
\]
and 
\[
\lambda_i(\alpha) = M(x_i(\alpha),\alpha)
\] 
where $x_1(\alpha)$ and $x_2(\alpha)$ are the unique roots of
Equation (\ref{ux}) in $[0,1]$ and in $[1,1/\alpha]$ respectively, then
we have the following lemma:
\begin{lem} \label{decreasing ga}
For a fixed $g\in [0,1]$, the function $\lambda_1$ is an
increasing function and the function $\lambda_2$ is an
decreasing function over $[0,\frac{1}{2}]$.
\end{lem}
\begin{proof}
  First, let us note that the function $f$ being strictly increasing
  over $[0,1]$ and strictly decreasing over $[1,1/\alpha]$, the implicit
  function theorem can be used to define the implicit function $\alpha
  \mapsto x_i(\alpha) \in [0,1]$ given by the equation
  $g=f(x)=x^{\alpha}\left( \frac{1-\alpha x}{1-\alpha
  }\right)^{1-\alpha}$, where $g$ is fixed. These functions are
  differentiable, and their derivative can be computed by implicitly
  differentiating the equation. In fact, taking the natural logarithm of
  Equation (\ref{ux}) and the derivative with respect to $\alpha$, after
  suitable simplifications, we obtain
\[
x'(\alpha) = \frac{x}{\alpha} \cdot \left(-1-\frac{1-\alpha x}{1-x}
\cdot \ln \left( 1-\frac{1-x}{1-\alpha x}\right) \right).
\]
Using the chain rule, we have
\[
\lambda'(\alpha) = \frac{\partial M}{\partial x}(x(\alpha),\alpha) \cdot
x'(\alpha) + \frac{\partial M}{\partial \alpha}(x(\alpha),\alpha).
\]
After suitable simplifications, we obtain
\begin{eqnarray*}
\frac{\partial M}{\partial \alpha}(x,\alpha) & = & -M^2(x,\alpha)
\cdot \frac{(1-x)^2}{x(1-\alpha x)^2}\\ \frac{\partial M}{\partial
  x}(x,\alpha) & = & -M^2(x,\alpha) \cdot\frac{-\alpha (1-x)
  (1-(2\alpha-1)x))}{x^2(1-\alpha x)^2}
\end{eqnarray*}
which leads to
\[
\lambda'(\alpha) =-M^2(x,\alpha) \cdot \frac{1-(2\alpha-1)x}{x(1-\alpha x)}
\left(\frac{1-x}{1-\alpha x}+
 \ln \left( 1-\frac{1-x}{1-\alpha x}\right) +
 \frac{(1-x)^2}{(1-\alpha x)(1-(2\alpha-1)x)} \right).
\] 
Note that since $\alpha\leq 1/2$, $\frac{1-(2\alpha-1)x}{x(1-\alpha
x)}\geq 0$. An application of Lemma \ref{gamma} shows that
 $\lambda_1'\geq 0$ and $\lambda_0'\leq 0$. This finishes the proof of
 the lemma.
\end{proof}

Let us now focus on the case where both $a$ and $h$ are known, and
an upper and a lower bound of $g$ is to be found, when $n=2$. The
problem can now be formulated as follows. Given two
real numbers $0<\alpha \leq 1/2$ and $0<h<1$, we want to find the
minimal and the maximal value of
\[
G(x,y) = x^{\alpha} y^{1-\alpha}
\]
where $x$ and $y$ fulfill the conditions
\begin{equation*}
\alpha x + (1-\alpha) y = 1 \; \mbox{ and } \;
\left( \frac{\alpha}{x} + \frac{1-\alpha}{y} \right)^{-1} =  h.
\end{equation*}
These two conditions imply that
\begin{equation} \label{condition3}
h=\left( \frac{\alpha}{x} + \frac{(1-\alpha)^2}{1-\alpha
  x} \right)^{-1}=\frac{x(1-\alpha x)}{\alpha-2\alpha x +x} = \varphi^2(x)
\end{equation}
and
\begin{equation*}
G(x,y)= f(x) = x^{\alpha} \left( \frac{1-\alpha x}{1-\alpha
}\right)^{1-\alpha}.
\end{equation*}
If $x_1$ and $x_2$ are the two unique solutions of Equation
(\ref{condition3}) with $x_1<1<x_2$, Lemma \ref{f-g} implies that
$f(x_1) > \varphi(x_1) = \varphi(x_2) > f(x_2)$ because of the
respective positions of $f$ and $\varphi$. This directly gives the
following lemma:
\begin{lem} \label{2-case ha}
Let $0<\alpha\leq1/2$ and $0<h<1$. If $x$ and $y$ fulfill the
conditions
\begin{equation*}
\alpha x + (1-\alpha) y = 1 \; \mbox{ and } \;
\left( \frac{\alpha}{x} + \frac{1-\alpha}{y} \right)^{-1} =  h,
\end{equation*}
then
\[
x_2^{\alpha} \left( \frac{1-\alpha x_2}{1-\alpha
}\right)^{1-\alpha}
\leq x^{\alpha} y^{1-\alpha} \leq
x_1^{\alpha} \left( \frac{1-\alpha x_1}{1-\alpha
}\right)^{1-\alpha}
\]
where $x_1$ and $x_2$ are the unique solutions over $[0,1]$ and
$[1,1/\alpha]$ respectively of the equation
\begin{equation}\label{vx}
h=\frac{x(1-\alpha x)}{\alpha-2\alpha x +x}.
\end{equation}
\end{lem}
As before, we would like now to prove that for a fixed $h$, and as a
function of $\alpha$, the minimal and maximal values above $G(x_2)$
and $G(x_1)$ are increasing and decreasing functions
respectively. More precisely, if we set
\[
N(x,\alpha) = x^{\alpha} \left( \frac{1-\alpha x}{1-\alpha
}\right)^{1-\alpha}
\]
and 
\[
\gamma_i(\alpha) =N(x_i(\alpha),\alpha)
\]
where $x_1(\alpha)$ and $x_2(\alpha)$ are the uniques root of Equation
(\ref{vx}) in $[0,1]$ and $[1,1/\alpha]$ respectively, then we have the
following lemma:
\begin{lem} \label{decreasing ha}
For a fixed $h\in [0,1]$, the function $\gamma_1$ is a
decreasing function and the function $\gamma_2$ is an
increasing function over $[0,1/2]$.
\end{lem}
\begin{proof}
The same argument used in the proof of Lemma \ref{decreasing ga} (with
$\varphi$ instead of $f$) shows that the function $x(\alpha)$ is well
defined, and after suitable simplifications, has the following derivative
\[
x'(\alpha) = \frac{x}{\alpha} \cdot \frac{1-x}{(1-2\alpha)x+1}.
\]
Using the chain rule, we have
\[
\gamma'(\alpha) = \frac{\partial N}{\partial x}(x(\alpha),\alpha) \cdot
x'(\alpha) + \frac{\partial N}{\partial \alpha}(x(\alpha),\alpha).
\]
After suitable simplifications, we obtain
\begin{eqnarray*}
\frac{\partial N}{\partial \alpha}(x,\alpha) & = & N(x,\alpha) \cdot
\left(\ln\left( 1- \frac{1-x}{1-\alpha x} \right) +
\frac{1-x}{1-\alpha x}  
 \right)\\
\frac{\partial N}{\partial x}(x,\alpha) & =  & N(x,\alpha) \cdot
\frac{\alpha}{x} \cdot \frac{1-x}{1-\alpha x} 
\end{eqnarray*}
which leads to
\[
\gamma'(\alpha) = N(x,\alpha) \cdot \left(  \frac{1-x}{1-\alpha x}+
 \ln \left( 1-\frac{1-x}{1-\alpha x}\right) +
\frac{(1-x)^2}{(1-\alpha x)(1-(2\alpha-1)x)}\right).
\] 
A straightforward application of Lemma \ref{gamma} shows that
$\gamma_1'\leq 0$ and $\gamma_2'\geq 0$ which finishes the proof of
the lemma.
\end{proof}

\Section{Proof of Theorems \ref{GA} and \ref{HA}}

The proofs of the two theorems are similar, so we treat them as a
whole and make the differences precise when needed. Without loss of
generality, we can suppose $n\geq 3$. In Theorem \ref{GA} (resp.
Theorem \ref{HA}), we suppose that the geometric mean $g>0$
(resp. harmonic mean $h>0$) and the arithmetic mean $a>0$ of a list of
$n$ strictly positive reals are given and we want to find sharp bounds
on the harmonic mean (resp. geometric mean). Before going further, let
us notice that the expression $\left( \sum_i \frac{\alpha_i}{x_i}
\right)^{-1}$, defined for $x_i>0$, can be continuously continued on
$\lbrack 0, \infty \lbrack^n$ by setting its value to $0$ as soon as
$x_i=0$ for some $i$. Let $\R_{\geq 0} = \lbrack 0, \infty \lbrack$
and let us define the three sets $S_h,S_g$ and $S_a$ as follows:
\[
S_a = \left\{ x\in \R^n \; | \; x_i \in \lbrack 0,
1/\alpha_i \rbrack , \, \sum \alpha_i x_i = 1
\right\},
\]
and
\[
S_g = \left\{x\in \R_{\geq 0}^n\; | \; \prod x_i^{\alpha_i} = g
\right\},  \; S_h = \left\{x\in \R_{\geq 0}^n\; | \;  \left( \sum
\frac{\alpha_i}{x_i}\right)^{-1} = h \right\}.
\]
The condition sets $C_1$ and $C_2$ on the $x_i$ related to Theorems
\ref{GA} and \ref{HA} respectively are given by $C_1 = S_g \cap S_a$
and $C_2 = S_h \cap S_a$. Because they are defined through the
preimage of closed sets via continuous maps, the sets $S_g$ and $S_h$
are closed and $S_a$ is compact. Therefore $C_1$ and $C_2$ are compact
in $\R^n$ as the intersection between a compact and a closed
set. Since the functions to optimize are well defined and continuous
on these sets, their maximum and minimum are reached, and we will
explicitly find them. The constraints being of class $C^1$, we use the
Lagrange multipliers to find these optimums. When the expression to
optimize is $\sum_{i=1}^n \frac{\alpha_i}{x_i}$ and the geometric and
the arithmetic means are known, the Lagrange's conditions are
\[
\frac{\partial}{\partial x_i} \left(\sum_{i=1}^n
\frac{\alpha_i}{x_i} -A \cdot \left( \sum_{i=1}^n \alpha_i x_i -1
\right)- B \cdot\left( \sum_{i=1}^n \alpha_i \ln(x_i) -\ln(g) \right)
\right) = 0
\]
which gives
\[
-\frac{1}{x_i^2}-A-\frac{B}{x_i}=0, \; \forall i=1,...,n.
\] 
When the expression to optimize is $\prod_{i=1}^n x_i^{\alpha_i}$
and the harmonic and the arithmetic means are known, the Lagrange's
conditions applied to the natural logarithm of the product are
\[
\frac{\partial}{\partial x_i} \left( \sum_{i=1}^n \alpha_i \ln(x_i)
-A \cdot \left( \sum_{i=1}^n \alpha_i x_i -1 \right)- B \cdot\left(
\sum_{i=1}^n \frac{\alpha_i}{x_i} - h^{-1}  \right) \right) = 0
\]
which gives
\[
\frac{1}{x_i}-A-\frac{B}{x_i^2}=0, \; \forall i=1,...,n.
\] 
In both cases, each $x_i$ is equal to one of the roots, say $X,Y$,
of a second degree polynomial. Since we supposed that the $x_i$'s are
not all equal, we have $X\neq Y$. Now, note that if $\alpha =
\sum_{i\in I} \alpha_i$ with $I=\{i | x_i = X\}$, then $1-\alpha =
\sum_{j\in J} \alpha_j$ with $J=\{j | x_j = Y\}$, and we may suppose
without loss of generality that $\alpha \in [0,1/2]$. Moreover
\begin{eqnarray*}
\sum_{i=1}^n \alpha_i x_i & = & \alpha X + (1-\alpha) Y = 1\\
\prod_{i=1}^n x_i^{\alpha_i} & = & X^{\alpha} Y^{1-\alpha}  = g\\
\left( \sum_{i=1}^n \frac{\alpha_i}{x_i} \right)^{-1} & = &
\left(  \frac{\alpha}{X}+ \frac{1-\alpha}{Y} \right)^{-1} = h.
\end{eqnarray*}
Suppose both the geometric and the arithmetic means are known and the
minimum and the maximum of the harmonic mean have to be
determined. Making use of Lemma \ref{2-case ga} and \ref{decreasing
  ga} and the previous notations, since $\alpha\leq 1/2$, we have
$H(X)< h < H(Y)$ where $X<1<Y$. The functions $H(X)$ and $H(Y)$ being
decreasing and increasing functions of $\alpha$, the minimum of $H(X)$
and the maximum of $H(Y)$ are reached when
$\alpha=\min_i\{\alpha_i\}$, $\alpha=0$ being impossible.\\

\noindent Similarly, suppose both the harmonic and the arithmetic
means are known and the minimum and the maximum of the geometric mean
have to be determined. Making use now of Lemma \ref{2-case ha} and
\ref{decreasing ha}, since $\alpha\leq 1/2$, we have $G(Y)< h < G(X)$
where $X<1<Y$. The functions $G(Y)$ and $G(X)$ being decreasing and
increasing functions of $\alpha$, the minimum of $G(Y)$ and the
maximum of $G(X)$ are reached as before when
$\alpha=\min_i\{\alpha_i\}$.\\

\noindent The statement of each Theorem \ref{GA} and \ref{HA} follows
then directly from the statements of Lemma \ref{2-case ga} and
\ref{2-case ha}.

\Section{Applications}

\begin{exmp}\label{ex3}
The first application is a bound on the trace of the inverse of a
matrix whose eigenvalues are all positive. This problem has been
treated by several authors (see \cite{bai} and the reference
therein). If $\lambda_i$ are the eigenvalues of such an $n \times n$
matrix $A$, then $\det(A) = \prod_{i=1}^n \lambda_i$, $\trace(A) =
\sum_{i=1}^n \lambda_i$, and $\trace(A^{-1}) = \sum_{i=1}^n
1/\lambda_i$. The connection with the arithmetic, geometric and
harmonic means is clear, and Corollary \ref{co} shows that
\[
\trace(A^{-1}) < e \cdot \left( \frac{\trace(A)}{n} \right)^n \cdot
\frac{1}{\det(A)} +\frac{n^2}{\trace(A)}.
\]
\end{exmp}

\begin{exmp}\label{ex4}
The second application is a bound on the quotient of some coefficients
of polynomials with positive roots. It has been known since Frans\'en
and Lohne \cite{fransen} (see also \cite{ineq}) that if the polynomial
\[
a_0x^n+a_1x^{n-1}+...+a_{n-1}x+a_{n} 
\]
has positive roots, then
\[
|a_{n-1}|\geq n^2 \left| \frac{a_0a_n}{a_1} \right|.
\]
An application of Corollary \ref{co} shows that the following reverse
inequality holds:
\[
|a_{n-1}| \leq n^2 \left| \frac{a_0a_n}{a_1} \right| + e |a_0|
\left|\frac{a_{1}}{na_0} \right|^n.
\] 
Indeed, if $\lambda_i$ are the roots of the polynomial, then
$|a_n/a_0| = \prod_{i=1}^n \lambda_i$, $|a_1/a_0| = \sum_{i=1}^n
\lambda_i$, and $|a_{n-1}/a_n| = \sum_{i=1}^n 1/\lambda_i$.
\end{exmp}

\vspace{3mm}

\noindent{\bf Acknowledgments.} The authors would like to thank the
reviewer for several helpful suggestions.

\end{document}